\documentclass[12pt]{article}

\usepackage{aeguill}
\usepackage[francais,english]{babel}
\usepackage[utf8]{inputenc}
\usepackage[T1]{fontenc}
\usepackage{amsfonts}
\usepackage{amssymb}
\usepackage{xspace}
\usepackage{graphicx}
\usepackage{amsmath}
\usepackage{amsthm}
\usepackage{stmaryrd}
\usepackage{MnSymbol}
\usepackage{endnotes}
\usepackage{yfonts}
\usepackage{makeidx}
\usepackage{multind}
\usepackage{fancyhdr}
\usepackage{pgf}
\usepackage{tikz}
\usepackage{dsfont}
\usepackage{bbm}
\usepackage[a4paper]{geometry}

\usetikzlibrary{positioning,arrows}

\theoremstyle{plain}
\newtheorem{thm}{Th\'eor\'eme}[section]
\newtheorem{prop}[thm]{Proposition}
\newtheorem{lem}[thm]{Lemme}

\theoremstyle{definition}
\newtheorem{df}[thm]{D\'efinition}
\newtheorem{dfp}[thm]{D\'efinition/Proposition}

\newenvironment{dem}{\paragraph{Preuve}}{\qed\bigskip}

\newcommand{\N}{\mathbb{N}}
\newcommand{\Z}{\mathbb{Z}}

\newcommand{\C}{\mathbb{C}}
\newcommand{\R}{\mathbb{R}}

\newcommand{\dd}[1]{\ensuremath{\operatorname{d}\!{#1}}}
\newcommand{\sui}{_{n\geq 1}}
\newcommand{\fcar}{\mathbbm{1}}
\newcommand{\n[1]}{||#1||}
\newcommand{\tend}{\rightarrow}

\newcommand{\ab[1]}{|#1|}
\newcommand{\app[5]}{#1\: \left\{\begin{array}{ccc} #2 & \longrightarrow & #3 \\ #4 & \longmapsto & #5 \end{array}\right.}

\newcommand{\ninf[1]}{\n[#1]_{\infty}}
\newcommand{\nti}{n\tend\infty}

\newcommand{\Ec}{E_{\C}}
\newcommand{\Hg}{\mathcal{H}^{\gamma}(\flag)}
\newcommand{\pth}{\phi_{\theta}}
\newcommand{\pthi}{\phi_{i\theta}}
\newcommand{\lth}{\lambda_{\theta}}
\newcommand{\lthi}{\lambda_{i\theta}}
\newcommand{\optr}{P_{\theta}}
\newcommand{\optri}{P_{i\theta}}
\newcommand{\mukx}{\mu_{k,\,\eta}}
\newcommand{\Ae}{A_e}

\newcommand{\Pj[1]}{\mathbb{P}(#1)}

\newcommand{\ym}{\Gamma_{\mu}}

\newcommand{\ma}{marche al\'eatoire\xspace}
\newcommand{\maX}{marche al\'eatoire sur $X$ associ\'ee \`a $G$ et $\mu$\xspace}
\newcommand{\Lie[1]}{\mathfrak{#1}}
\newcommand{\flag}{\mathcal{P}}
\newcommand{\muen}{\mu^{*n}}

\newcommand{\sm}{\sigma_{\mu}}
\newcommand{\smb}{\overline{\sigma_{\mu}}}
\newcommand{\tcl}{th\'eor\`eme central-limite\xspace}
\newcommand{\cad}{c\textquoteright est-\`a-dire\xspace}
\newcommand{\Bf}{B\times \flag}
\newcommand{\stb}{\tilde{\sigma}_k}
\newcommand{\sgb}{\overline{\sigma}}

\newcommand{\stX}{sur tout $X$\xspace}
\newcommand{\Bb}{\tilde{B}}

\newcommand{\mbb}{\tilde{\beta}}

\newcommand{\borBf}{\mathcal{B}^{\flag}}

\newcommand{\borBb}{\tilde{\mathcal{B}}}
\newcommand{\Lun}{\mathrm{L}^1}
\newcommand{\Linf}{\mathrm{L}^{\infty}}
\newcommand{\nLun[1]}{||#1||_{\Lun}}

\DeclareMathOperator{\Ker}{Ker}
\DeclareMathOperator{\Ima}{Im}

\DeclareMathOperator{\Supp}{Supp}
\DeclareMathOperator{\Stab}{Stab}
\DeclareMathOperator{\End}{End}

\begin{document}

\title{Un critère de récurrence pour certains espaces homogènes}
\author{Caroline Bruère}
\date{Mai 2016}

\maketitle

\selectlanguage{english}

\begin{abstract} Let $G$ be a real connected algebraic semi-simple Lie group, 
and $H$ an algebraic subgroup of $G$.
Let $\mu$ be a probability measure on $G$, with finite exponential moment,
whose support spans a Zariski-dense subsemigroup of $G$.
Let $X=G/H$ be the quotient of $G$ by $H$. 
We study the Markov chain on $X$ with transition probability $P_x=\mu *\delta_x$ for $x\in X$. 
We prove that either for every $x\in X$, almost every trajectory starting from $x$ is transient 
or for every $x\in X$, almost every trajectory starting from $x$ is recurrent.
In fact, this recurrence is uniform over all $X$, 
i.e. there exists a compact set $C\subset X$ such that
for each point $x\in X$, every trajectory starting in $x$
almost surely returns to $C$ infinitely often.
Furthermore, we give a criterion for recurrence depending on $G$, $H$, and $\mu$.  
\end{abstract}

\selectlanguage{francais}

\begin{abstract} Soit $G$ un groupe de Lie algébrique connexe semi-simple réel,
$H$ un sous-groupe algébrique de $G$, 
$\mu$ une mesure de probabilité sur $G$ à moment exponentiel fini
dont le support engendre un sous-semi-groupe Zariski-dense de $G$.
Soit $X=G/H$ le quotient de $G$ par $H$.
On étudie la chaîne de Markov sur $X$ de probabilité de transition $P_x=\mu *\delta_x$ pour $x\in X$.
On montre que soit pour tout $x\in X$, presque toute trajectoire partant de $x$ est transiente,
soit pour tout $x\in X$, presque toute trajectoire partant de $x$ est récurrente.
Cette récurrence est en fait uniforme, 
c'est-à-dire que pour tout point $x\in X$, 
presque toute trajectoire partant de $x$ revient infiniment souvent 
dans un compact $C\subset X$ ne dépendant pas de $x$. 
De plus, on donne un critère de récurrence en fonction de $G$, $H$, et $\mu$.

\end{abstract}

\pagebreak

\tableofcontents

\section{Introduction}

Un célèbre théorème de Pólya (1928) dit qu'étant donnés $d\geq 1$ un entier 
et $\mu$ une mesure de probabilité sur $\R^d$, 
à moment exponentiel fini, 
centrée, 
la \ma correspondante sur $\R^d$ est récurrente, 
si et seulement si $d$ vaut $1$ ou $2$.
Nous allons démontrer un résultat analogue dans le cas d'une marche aléatoire 
sur certains espaces homogènes associés à des groupes de Lie semi-simples.

Soit $G$ un groupe topologique localement compact à base dénombrable, 
$H$ un sous-groupe fermé de $G$, 
et $X=G/H$ le quotient de $G$ par $H$.
Soit $\mu$ une mesure de probabilité sur $G$.
La \emph{\maX} est la chaîne de Markov sur l'espace $X$ de probabilité de transition 
$P_x=\mu *\delta_x$ pour $x\in X$. 
Notons $B=G^{\N^*}$, et $\beta=\mu^{\N^*}$ la mesure de probabilité produit sur $B$. 

\begin{df} \label{defrectrans}
On dit que la \ma sur $X$ est \emph{récurrente} en un point $x\in X$ 
s'il existe un compact $C$ de $X$ tel que 
\[\beta(\{b\in B\, |\, \forall n_0\in\N, \exists n\geq n_0 : b_n\cdots b_1x\in C\})=1.\] 
On dit que la marche est \emph{transiente} en un point $x\in X$ 
si pour tout compact $C$ de $X$, on a
\[\beta(\{b\in B\, |\, \exists n_0\in\N, \forall n\geq n_0 : b_n\cdots b_1x\notin C\})=1.\]
On dit que la \ma est récurrente (respectivement transiente) sur tout $X$ si elle l'est en tout point. 
On dit que la \ma sur $X$ est \emph{uniformément récurrente sur tout $X$} 
s'il existe un compact $C$ de $X$ tel que pour tout $x\in X$
\[\beta(\{b\in B\, |\, \forall n_0\in\N, \exists n\geq n_0 : b_n\cdots b_1x\in C\})=1.\] 
\end{df}

Remarquons qu'avec ces définitions, 
un point $x\in X$ peut être récurrent 
sans que presque toute trajectoire revienne dans tout voisinage de $x$. 
Par exemple, dans le cas de l'action sur l'espace projectif 
de dimension $m-1$ d'un sous-groupe de Schottky Zariski-dense de $\mathrm{SL}(m)$, 
engendré par deux matrices,
l'espace $\omega$-limite de presque toutes les trajectoires est un compact de Cantor ; 
tous les points de l'espace projectif sont récurrents selon notre définition.

Dans le cas où $\mu$ est étalée, c'est-à-dire absolument continue par rapport à la mesure de Haar, 
il existe des théorèmes de dichotomie, 
i.e. des conditions pour que les états soient tous récurrents ou tous transients,
notamment le théorème de Hennion et Roynette (\cite{HennionRoynette}),
 affiné par Elie dans \cite{ElieHR}, 
et dont une preuve plus courte est fournie par Revuz dans \cite{RevuzHR}. 
Ce théorème concerne une classe très large d'espaces homogènes.
Cependant, la condition "$\mu$ étalée" est très restrictive.
Elle ne permet par exemple pas de traiter le cas décrit ci-dessus 
d'une mesure à support fini engendrant un semi-groupe discret Zariski-dense dans $G$. 
Nous allons donner un critère de récurrence n'utilisant pas cette condition, 
sur une classe plus restreinte d'espaces homogènes. 

Considérons pour $G$ un groupe de Lie algébrique connexe semi-simple réel, 
et pour $H$ un sous-groupe algébrique de $G$. 
Munissons $G$ d'une mesure de probabilité $\mu$, 
dont le support engendre un semi-groupe $\ym$ Zariski-dense dans $G$. 
Mentionnons tout d'abord les travaux sur la récurrence sur les espaces homogènes 
de Guivarc'h et Raja (\cite{GuiRaj11} et \cite{GuiRaj12}), 
de Benoist et Quint (\cite{BenQuiFinVol},  \cite{BenQuiIntro}, et \cite{BenRecLat}).
Plusieurs problèmes apparaissent. 
Les trajectoires issues de chaque point ont-elles presque toutes le même comportement,
i.e. la marche est-elle soit récurrente, soit transiente en chaque point ? 
Si c'est le cas, a-t-on un théorème de dichotomie, 
i.e. la marche est-elle soit récurrente sur tout $X$, soit transiente sur tout $X$? 
Enfin, quels critères permettent, selon $G$, $H$, et $\mu$, 
de dire si la \ma est transiente ou récurrente sur tout $X$ ? 

Nous allons répondre par l'affirmative aux deux premières questions, 
dans le cadre défini plus haut, pour des mesures à moment exponentiel fini. 
Cependant, dans une classe plus générale d'espaces homogènes, 
on peut construire un exemple où la marche est récurrente en certains points, 
et transiente en d'autres (voir \cite{Bru2}).
Nous donnerons également une condition nécessaire et suffisante de récurrence. 
Commençons par une définition :

\begin{df} Soit $G$ un groupe de Lie algébrique connexe semi-simple réel, 
soit $\mu$ une mesure de probabilité sur $G$. 
On dit qu'elle est \emph{à moment exponentiel fini} 
s'il existe $\alpha \in \R_+^*$ tel que le moment d'ordre exponentiel $\alpha$ soit fini :
\[ \int_G \n[g]^{\alpha}\dd \mu(g) < \infty .\]
On dit qu'elle est \emph{Zariski-dense} si son support engendre un sous-semi-groupe $\ym$ 
de $G$ qui est Zariski-dense dans $G$. 
\end{df}

Énonçons alors le théorème de dichotomie et le critère de récurrence.

\begin{thm}\label{thmdicho} 
\textbf{(Théorème de dichotomie pour des marches aléatoires sur certains espaces homogènes)}
Soit $G$ un groupe de Lie algébrique connexe semi-simple réel, 
et $H$ un sous-groupe algébrique de $G$. 
Notons $X$ l'espace homogène $X=G/H$.
Soit $\mu$ une mesure de probabilité sur $G$ à moment exponentiel fini, et Zariski-dense. 
Alors la \maX est soit uniformément récurrente sur tout $X$, soit transiente sur tout $X$. \end{thm}

Ce théorème est en fait un corollaire du théorème \ref{thmclass} suivant, 
qui décrit les cas où la \maX est récurrente, et ceux où elle est transiente sur tout $X$. 
Soit $N$ un sous-groupe unipotent maximal de $G$, 
$P$ le normalisateur de $N$, 
$A$ un sous-tore déployé maximal de $P$,
et $\Lie[a]$ l'algèbre de Lie de $A$.
Notons $\flag=G/P$ la variété drapeau de $G$, 
$\sigma : G\times \flag \rightarrow \Lie[a]$ le cocycle d'Iwasawa (cf. \ref{defcocycle}),
et $\sm$ sa moyenne, aussi appelée le vecteur de Lyapounov de $\mu$.

\begin{thm}\label{thmclass} 
\textbf{(Critère de récurrence pour certains espaces homogènes)}
Avec les hypothèses et notations ci-dessus,  
la \maX est uniformément récurrente \stX
si et seulement si $H$ contient un conjugué de $A'N$, 
où $A'$ est un sous-groupe de $A$ de codimension au plus $2$, 
dont l'algèbre de Lie $\Lie[a']$ contient la moyenne $\sm$ du cocycle d'Iwasawa.
Sinon, la marche est transiente \stX. 
 \end{thm}

La preuve des théorèmes sera donnée en \ref{preuveprinc}.
Remarquons d'abord qu'à conjugaison près, 
les sous-groupes algébriques $H$ de $G$ contiennent 
un sous-groupe cocompact de la forme $A'N'$, 
avec $A'$ et $N'$ des sous-groupes algébriques de $A$ et $N$ respectivement.
Nous commencerons par montrer -- c'est l'objet de la proposition \ref{roledeN} -- que, 
si $N'$ est différent de $N$, 
la \ma est nécessairement transiente \stX.
Cela permettra de se ramener à des groupes de la forme $A'N$,
et à la récurrence ou transience d'un cocycle, 
le cocycle d'Iwasawa, sur un espace vectoriel réel.
La démonstration utilisera des propriétés des opérateurs de transferts 
développées dans \cite{Livre},
et des arguments classiques de la démonstration 
du théorème de récurrence de Pólya (voir par exemple \cite{Spitz}).
Un théorème de Conze et Schmidt (\cite{Conz}, \cite{Schm98}),
qui donne la récurrence de certains cocycles 
si leur distribution asymptotique est normale,
permettra ensuite de montrer une condition suffisante de récurrence,
grâce au \tcl pour le cocycle d'Iwasawa (\cite{Goldguiv}).
On montrera enfin que la récurrence est uniforme
en montrant l'ergodicité d'un système dynamique fibré,
grâce à des travaux de Guivarc'h sur l'exactitude (\cite{Guiv89}).

\section{Définitions et préliminaires}

\subsection{La décomposition d'Iwasawa}

Prenons à nouveau
$G$ un groupe de Lie algébrique connexe semi-simple réel.
Notons $A$ un tore déployé maximal de $G$,
$N$ un sous-groupe unipotent maximal de $G$, 
dont le normalisateur $P$ contient $A$, et
$\Ae$ la composante connexe de $A$ contenant l'identité.
Notons $\Lie[g]$ l'algèbre de Lie associée au groupe $G$,
$\Lie[a]$ celle associée au tore $A$,
et $\Sigma$ l'ensemble des racines restreintes de $\Lie[g]$, sous l'action adjointe de $\Lie[a]$.
Notons $\Sigma_+$ l'ensemble de racines positives
associé au choix de $A$ et $N$,
$\Lie[a]_+\subset \Lie[a]$ la chambre de Weyl associée à $\Sigma_+$. 
Notons également $A_+=\exp \Lie[a]_+\subset A$.

\begin{dfp} Il existe une écriture de $G$ sous la forme 
\[G=K\Ae N\]
où $K$ est un sous-groupe compact de $G$,
appelée la \emph{décomposition d'Iwasawa} de $G$.
$G$ s'écrit alors également 
\[G=KA_+ K,\]
sa \emph{décomposition de Cartan}.
\end{dfp}

On pourra voir \cite[ch. 5]{Livre} pour une preuve de l'existence de ces décompositions.
La proposition suivante donne une décomposition analogue 
à la décomposition d'Iwasawa
pour tout sous-groupe algébrique $H$ de $G$.

\begin{prop}\label{iwasH} Soit $H$ un sous-groupe algébrique de $G$.
Il existe un sous-tore $A'$ de $A$ et un sous-groupe algébrique $N'$ de $N$
tels que, à conjugaison près, 
$A'N'$ soit un sous-groupe cocompact de $H$. 
\end{prop}

\begin{dem} $H$ est un sous-groupe de Lie algébrique de $G$ ; 
il a donc une décomposition de la forme 
\[H=K'A''N'',\]
avec $K'$ un sous-groupe compact de $H$, $A''$ un tore déployé, 
et $N''$ un sous-groupe unipotent maximal normalisé par $A''$. 
Puisque tous les sous-groupes trigonalisables maximaux d'un groupe réductif sont conjugués 
(\cite[Thm 8.2]{BorTit}),
$A''N''$ est conjugué à un sous-groupe de $AN$. 
Supposons, pour simplifier, $A''N''\subset AN$.
On a alors $N''\subset N$. 
D'après \cite[Thm 11.6]{BorTit}, tous les tores déployés maximaux de $AN$ 
sont conjugués dans $AN$.
$A''$ est donc conjugué à un sous-tore $A'$ de $A$ dans $AN$; 
comme $N$ est normalisé par $AN$, 
les conjugués de $N''$ par des éléments de $AN$ sont des sous-groupes de $N$. 
Il existe donc $N'\subset N$ tel que 
$A''N''$ est conjugué à $A'N'$. 
\end{dem}

Cette décomposition va permettre de définir un cocycle crucial 
pour l'étude des marches aléatoires sur l'espace homogène $X=G/H$.

\subsection{Le cocycle d'Iwasawa}

Notons $\flag=G/P$ la variété drapeau de $G$,
et notons comme auparavant $\mu$ une mesure de probabilité
sur $G$ Zariski-dense, à moment exponentiel fini.
Le fait suivant est dû à Furstenberg (\cite{Furst63}), 
Guivarc'h, Raugi (\cite{GuiRaug85}), Goldsheid et Margulis (\cite{GoldMarg}).

\begin{prop} Avec les hypothèses et notations ci-dessus, 
il existe sur $\flag$ une unique mesure $\mu$-stationnaire,
qu'on notera dorénavant $\nu$,
appelée \emph{mesure de Furstenberg}.
Rappelons qu'une mesure $\nu$ est dite $\mu$-stationnaire,
$\mu$-invariante, ou $\mu$-harmonique si elle vérifie la condition
\[\nu=\mu * \nu = \int_G g_*\nu\dd\mu(g).\]
\end{prop}

On peut alors définir le cocycle d'Iwasawa,
en notant $\Lie[a]$ l'algèbre de Lie de $A$ :

\begin{df}\label{defcocycle} Le \emph{cocycle d'Iwasawa} de $G$ est l'application  
$\sigma:G\times \flag \rightarrow \Lie[a]$ telle que pour tout $g\in G$, et pour tout $\eta\in\flag$, 
en notant $k$ un élément de $K$ tel que $\eta=kP$, on ait 
\[ gk \in K \text{exp}(\sigma(g,\eta)) N. \]
\end{df}

\begin{lem} Le cocycle d'Iwasawa est défini de manière unique, 
et c'est bien un cocycle : on a pour tous $g,\,g'\in G$, pour tout $\eta \in \flag$, l'égalité 
\[ \sigma(gg',\eta)=\sigma(g,g'\eta)+\sigma(g',\eta)  \]
\end{lem}

On trouvera une démonstration du lemme dans \cite{Livre}. 

\begin{df} On note alors $\sm$ la moyenne du cocycle d'Iwasawa: 
\[\sm=\int_{G\times\flag} \sigma(g,\,\eta)\dd \mu(g) \dd \nu(\eta). \]
On appelle également \emph{vecteur de Lyapounov de $\mu$} cette quantité $\sm$. 
\end{df}

\subsection{La fonction de Green}\label{utilisationcocycle}

On va fournir dans cette partie une condition suffisante de transience dans un cadre large.
Supposons plus généralement que  $G$ est un groupe topologique localement compact,
et $X$ un espace topologique localement compact muni d'une action continue de $G$.
Soit $\mu$ une mesure de probabilité borélienne sur $G$.
Soit $C$ un compact de $X$. 

\begin{df}
Notons $G_x^n(C)$ l'espérance du nombre de fois 
où une trajectoire $b_n\cdots b_1 x$ revient dans $C$ jusqu'au temps $n$ :
\[G_x^n(C)= \int_{b\in B} \sum_{k=0}^n \fcar_C(b_k\cdots b_1 x)\dd\beta(b).\]
La \emph{fonction de Green en $x\in X$ de la \maX } est l'application $G_x$ 
qui à un compact $C$ de $X$ associe 
l'espérance $G_x(C)$ du nombre de retours dans $C$ d'une trajectoire issue de $x$ :
\[ G_x(C)=\lim_{\nti} G_x^n(C) \in [0,\,\infty ] .\]
\end{df}

On peut en déduire la condition suffisante suivante de transience de la \maX en un point $x\in X$: 

\begin{lem}\label{Gxtrans} Soit $x\in X$. Si pour tout compact $C$ de $X$, 
la fonction de Green $G_x(C)$ est finie, i.e. 
\[ \lim_{\nti} G_x^n(C) < \infty, \]
alors $x$ est transient. 
\end{lem}

\begin{dem} Soit $x\in X$ tel que pour tout compact $C$ de $X$, 
la fonction de Green $G_x(C)$ soit finie.
Soit $C$ un compact de $X$. 
Ecrivons la quantité $G_x(C)$ : 
\begin{align*}
G_x(C) &= \sum_{k=0}^\infty \int_{b\in B}\fcar_C(b_k\cdots b_1 x) \dd\beta (b)\\
& = \int_{b\in B} \sum_{k=0}^\infty  \fcar_C(b_k\cdots b_1 x) \dd\beta (b) \\
\end{align*}
Comme $G_x(C)$ est finie, on en déduit que pour $\beta$-presque tout $b\in B$, 
la somme $\sum_{k=0}^\infty  \fcar_C(b_k\cdots b_1 x)$ est finie, et $x$ est bien transient. 
\end{dem}

\section{Une première condition suffisante de transience}

On suppose de nouveau que $G$ est un groupe de Lie algébrique connexe semi-simple réel, 
et que $H$ est un sous-groupe algébrique de $G$.

\subsection{Un critère de transience}

La proposition suivante va permettre de simplifier le problème, 
et de se ramener au cas où $H$ contient $N$. 

\begin{prop}\label{roledeN} Si le sous-groupe $H$ de $G$ ne contient pas de conjugué de $N$, 
alors la \ma sur $X=G/H$ associée à $G$ et $\mu$ est transiente sur tout $X$.
\end{prop}

\paragraph{Preuve de la proposition \ref{roledeN}} Raisonnons par l'absurde 
en supposant que $H$ ne contient pas $N$, et que la \maX n'est pas transiente. 
Rappelons le théorème de Chevalley sur une représentation du quotient $G/H$, 
dont on pourra trouver une preuve par exemple dans \cite[pg 38]{BenResGL}.

\begin{prop} (Chevalley)\label{repchev} Soit $G$ un groupe algébrique réel 
et $H\subset G$ un sous-groupe algébrique de $G$. 
Alors il existe une représentation algébrique de $G$ dans un $\R$-espace vectoriel $V$ 
et un point $x_0\in\Pj[V]$ tel que le stabilisateur $\Stab x_0$ de $x_0$ est égal à $H$. 
\end{prop}

Notons $X_0$ l'ensemble des points de $X$ invariants par un conjugué de $N$. 
L'ensemble $X_0$ est soit $X$ tout entier, soit l'ensemble vide. 
On se ramène donc à montrer que $X_0$ contient au moins un élément.
Soit $(V,\,\rho)$ une représentation (réelle) de $G$ telle qu'il existe $x_0\in\Pj[V]$ 
vérifiant $H=\Stab x_0$ et telle que l'orbite $G\cdot x_0$ de $x_0$
sous l'action de $G$ engendre l'espace vectoriel $V$, donnée par la proposition \ref{repchev}.
On a alors l'isomorphisme de $G$-espaces $G\cdot x_0\simeq X$.
Notons $\Pj[V]^N$ l'ensemble des points de $\Pj[V]$ fixés par l'action de $N$. 
Alors on peut écrire $X_0\simeq G\cdot\Pj[V]^N\cap \, G\cdot x_0$.
Notons $B=G^{\otimes\N^*}$, et $\beta=\mu^{\otimes\N^*}$.
Pour $b\in B$, notons $b_n\cdots b_1=p_n(b)=p_n$, 
et $\pi$ une valeur d'adhérence de la suite $(\frac{\rho(p_n)}{\n[\rho(p_n)]})\sui$. 
Notons $(\frac{\rho(p_n)}{\n[\rho(p_n)]})_{n\in S}$ la sous-suite de 
$(\frac{\rho(p_n)}{\n[\rho(p_n)]})\sui$ qui converge vers $\pi$. 
On a besoin des deux lemmes suivants, qui sont démontrés dans \cite{Art} :

\begin{lem}\label{cvupn} La convergence vers $\pi$ de la suite 
$(\frac{\rho(p_n)}{\n[\rho(p_n)]})_{n\in S}$ est uniforme sur les compacts de 
$\Pj[V]\setminus\Pj[\Ker\pi]$. 
\end{lem}

\begin{dem} Ecrivons la décomposition de Cartan du groupe $G$ : 
\[G=K\Ae K.\] 

Les éléments $p_n$ peuvent s'écrire sous la forme 
\[p_n=k_{1n}a_n k_{2n},\]
où les $k_{i,\,n}$ sont dans $K$, et les $a_n$ sont dans $\Ae$. 
Puisque $K$ est compact, on peut supposer, 
quitte à prendre une sous-suite de $(p_n)_{n\in S}$, 
que les suites $(k_{in})_{n\in S}$ convergent vers des éléments $k_i$ de $K$.
Il existe alors un endomorphisme $a\in\End V$, tel que  
\[\lim_{\nti, \,n\in S}\frac{\rho(a_n)}{\n[\rho(a_n)]}=a\]
et 
\[\rho(\pi)=\rho(k_1)a\,\rho(k_2).\]
Ainsi, il suffit de montrer que la convergence de $(\rho(a_n))_{n\in S}$ vers $a$ 
est uniforme sur tout compact de $\Pj[V]\setminus\Pj[\Ker\pi]$ ; 
comme $A$ est un tore déployé, 
il existe une base orthonormée de $V$ dans laquelle les éléments de $\rho(A)$ sont diagonaux,
et on a la convergence uniforme sur tout compact voulue. 
\end{dem}

\begin{lem}\label{incim}Pour $\beta$-presque tout $b\in B$, 
pour toute valeur d'adhérence $\pi$ de la suite $(\frac{\rho(p_n)}{\n[\rho(p_n)]})\sui$, 
on a $\Pj[\Ima\pi]\subset G\cdot\Pj[V]^N$.
\end{lem}

\begin{dem} C'est en fait le \cite[Lemme 4.3]{Art}.
Rappelons qu'on note $\Lie[g]$ l'algèbre de Lie associée au groupe $G$,
et $\Sigma$ l'ensemble des racines restreintes de $\Lie[g]$, sous l'action adjointe de $\Lie[a]$.
Pour $\alpha\in\Sigma$, on note $\Lie[g]^{\alpha}$ l'espace radiciel associé dans $\Lie[g]$.
Notons $\Sigma_+$ l'ensemble de racines positives
associé au choix de $A$ et $N$,
et choisissons $\Pi\subset\Sigma_+$ un sous-ensemble de racines simples.
Notons $\Lie[n]=\oplus_{\alpha\in\Sigma_+}\Lie[g]^{\alpha}$ l'algèbre de Lie de $N$, 
et $\Lie[a]_+\subset \Lie[a]$ la chambre de Weyl associée à $\Sigma_+$. 
Choisissons une valeur d'adhérence $\pi$ de la suite $(\frac{\rho(p_n)}{\n[\rho(p_n)]})\sui$,
pour un $b\in B$.
Pour tout $n\in\N$,  étudions la décomposition de Cartan de $p_n$ : 
\[p_n=k_{1,\,n}\exp{X_n} k_{2,\,n},\]
avec $X_n\in\Lie[a]_+$, et  $k_{1,\,n},\,k_{2,\,n}\in K$.
Par compacité de $K$, on peut supposer que les suites $k_{1,\,n}$ et $k_{2,\,n}$
convergent vers des limites $k_{1,\,\infty},\, k_{2,\,\infty}$.
On se ramène donc à l'étude de la suite des $g_n=\exp{X_n}$. 
Notons $V=\oplus_{\chi\in\Sigma(\rho)} V_{\chi}$ la décomposition de $V$ en espaces de poids. 
On peut alors écrire $\rho(g_n)=(\exp{\chi(X_n)})_{\chi\in\Sigma(\rho)}$ pour $n\in\N$. 
Notons $(\frac{\rho(g_n)}{\n[\rho(g_n)]})_{n\in S}$ 
la sous-suite de $(\frac{\rho(g_n)}{\n[\rho(g_n)]})_{n\in \N}$ convergeant vers $\pi$.
Notons alors $\Sigma'\subset\Sigma(\rho)$ l'ensemble de caractères de $\Sigma(\rho)$ tels que
\begin{itemize}
\item pour tous $\chi,\,\chi'\in\Sigma'$, on ait 
\[\limsup_{ n\in S}\ab[\chi(X_n)-\chi'(X_n)]<\infty,\]
\item pour tous $\chi\in\Sigma'$, $\chi'\notin\Sigma'$, on ait
\[\lim_{n\in S}\chi(X_n)-\chi'(X_n)=+\infty.\]
\end{itemize}
En passant à la limite, on obtient l'inclusion
\[\Ima\pi\subset \oplus_{\chi\in\Sigma'}V_{\chi}.\]
Soient maintenant $\alpha\in\Sigma_+$, et $\chi\in\Sigma'$.
Considérons l'action de $\Lie[g]^{\alpha}$ sur $V_{\chi}$.
On a l'inclusion $\Lie[g]^{\alpha}\cdot V_{\chi}\subset V_{\alpha+\chi}$.
Mais, comme $\alpha$ est une racine positive,
presque tout $b\in B$ vérifie la propriété 
\[\lim_{\nti}\alpha(X_n)=\infty,\]
par la loi des grands nombres et le théorème de positivité
du premier exposant de Lyapounov (voir \cite[Thm 9.9]{Livre}).  
La racine $\chi+\alpha$
ne peut donc pas être un poids de $V$. 
On en déduit 
\[\Lie[g]^{\alpha}\cdot V_{\chi} = \{0\},\]
ce qui prouve le résultat. 
\end{dem}

Reprenons la démonstration de la proposition \ref{roledeN}.
Comme on suppose la \ma non-transiente, 
on peut choisir $x\in X$ tel qu'il existe un compact $C_0\subset X$
tel que la probabilité qu'une trajectoire issue de $x$ 
ne passe qu'un nombre fini de fois dans $C_0$ 
soit strictement plus petite que $1$, c'est-à-dire 
\[\beta(\{b\in B\, |\, \exists n_0\in\N, \forall n\geq n_0 : b_n\cdots b_1x\notin C_0\})<1.\]
Choisissons donc $b\in B$ tel que la trajectoire issue de $x$ 
associée revienne une infinité de fois dans $C_0$. 
La suite des $(b_n\cdots b_1 x)_{n\geq 1}$ a donc une valeur d'adhérence, 
$y\in C_0\subset X$ ; 
soit $(b_n\cdots b_1 x)_{n\in S}$ une sous-suite extraite convergeant vers $y$.
Extrayons une sous-suite de $(\frac{\rho(p_n)}{\n[\rho(p_n)]})_{n\in S}$ convergeant vers $\pi$, qu'on note $(\frac{\rho(p_n)}{\n[\rho(p_n)]})_{n\in S'}$.
Alors on obtient 
\[y=\lim_{\nti,\,n\in S'} \frac{\rho(p_n)(x)}{\n[\rho(p_n)]}=\pi(x).\]
D'après le lemme \ref{incim}, le point $y$ de $X\simeq G\cdot x_0$ 
est dans $G\cdot\Pj[V]^N\cap \, G\cdot x_0$, 
et donc $X_0$ est non-vide !
\qed

\bigskip

Il suffira donc d'étudier les cas où $H$ contient un conjugué de $N$.

\subsection{Simplification du problème}

Dans le cas où $H$ contient un conjugué de $N$, 
le lemme \ref{iwasH} fournit un sous-tore $A'$ de $A$ tel que, à conjugaison près, 
$A'N$ est un sous-groupe cocompact de $H$. 
La récurrence d'une \maX est alors liée au comportement du cocycle d'Iwasawa.
Notons $\Lie[a']$ l'algèbre de Lie de $A'$, et $E$ le quotient $ \Lie[a]/\Lie[a']$. 
Notons \[\sgb:G\times\flag\rightarrow E\] 
le cocycle $p\circ \sigma$, où $p:\Lie[a]\rightarrow E$ est la projection canonique. 
Notons 
\[\smb=p(\sm)=\int_{G\times\flag} \sgb(g,\,\eta)\dd \mu(g) \dd \nu(\eta) \] sa moyenne.
Notons $\pi$ la projection canonique $\pi:X\rightarrow \flag$. 
On obtient une action de $G$ sur $\flag\times E$ via la formule
\[g\cdot(\eta,\,t)=(g\cdot \eta, \, t+\sgb(g,\,\eta)),\]
pour tous $g\in G,\,\eta\in\flag,\,t\in E$.
En notant $\eta_0$ le point base de $\flag$, 
on obtient que $A'N$ est inclus dans le stabilisateur de $(\eta_0,\,0)$
sous l'action de $G$, et que l'application 
\[\app[]{G/A'N}{\flag\times E}{g}{g\cdot (\eta_0,\,0)}\]
est propre. 
Étudier la récurrence de la marche aléatoire induite par une mesure $\mu$
sur $G/H$ revient donc à étudier le cocycle $\sgb$ sur $\flag$. 
La récurrence et la transience de la \maX en un point $x\in X$ s'écrivent alors de la manière suivante :

\begin{lem} \label{reccocycleiwas} La \maX est récurrente en $x\in X$ si et seulement  
s'il existe un compact $C$ de $E$ tel que 
\[\beta(\{b\in B\, |\, \forall n_0\in\N, \exists n\geq n_0 : \sgb(b_n\cdots b_1, \,\pi(x))\in C\})=1.\] 
Elle est uniformément récurrente sur tout $X$ 
si et seulement s'il existe un réel $A>0$ tel que
\[\forall \eta\in \flag, \forall t\in E,\, 
\beta(\{b\in B\, |\, \forall n_0\in\N, 
\exists n\geq n_0 : \n[t+\sgb(b_n\cdots b_1, \,\eta)]\leq A\})=1.\]
Elle est transiente en $x\in X$ si et seulement si, pour tout compact $C$ de $E$, on a 
\[\beta(\{b\in B\, |\, \exists n_0\in\N, \forall n\geq n_0 : \sgb(b_n\cdots b_1,\,\pi(x))\notin C\})=1.\]
\end{lem}

En réécrivant le lemme \ref{Gxtrans}, on obtient une condition suffisante de transience de la \ma
utilisant le cocycle d'Iwasawa :

\begin{lem}\label{Gxtranscocycle}
Soit $x\in X$. 
Introduisons alors, pour tout compact de $C$ de $E$,  la quantité $F_x(C)$ :
\[F_x(C)=\sum_{k=0}^{\infty}\int_{g\in G}\fcar_{C}(\sgb(g,\,x))\dd\mu^{*k}(g)\]
Si pour tout compact $C$ de $E$, la fonction de Green modifiée $F_x(C)$ est finie :
\[F_x(C)<\infty\]
alors la \maX est transiente en $x$. 
\end{lem}

On va donc à présent s'intéresser à la récurrence ou transience du cocycle $\sgb$ sur $E$.

\section{Une seconde condition suffisante de transience}

On va supposer dans les deux prochaines parties que $H$ s'écrit sous la forme $H=A'N$. 
Rappelons qu'on note $\Lie[a]$ l'algèbre de Lie de $A$, 
$\Lie[a']$ l'algèbre de Lie de $A'$, 
et $E$ le quotient $\Lie[a]/\Lie[a']$. 
Notons $d$ la dimension de $E$, $\Ec$ son complexifié. 
Notons enfin à nouveau $\nu$ l'unique mesure $\mu$-invariante sur $\flag$. 
Considérons la projection du cocycle d'Iwasawa $\sgb: G\times \flag\rightarrow E$. 
Ce qui suit est dans l'esprit de la démonstration du théorème de récurrence de Pólya 
proposée dans \cite[\S 6, 7, 8]{Spitz}.

\subsection{Opérateur de transfert associé à une mesure}

Commençons par une définition.

\begin{df} On dit qu'une fonction $\phi : \flag\rightarrow \C$ est 
$\gamma$-Hölderienne pour $0<\gamma\leq 1$
si on a 
\[\sup_{x,\,x'\in\flag} \frac{\ab[\phi(x)-\phi(x')]}{d(x,\,x')^{\gamma}} <\infty.\]
On note $\Hg$ l'ensemble des fonctions $\gamma$-Hölderiennes sur $\flag$. 
L'\emph{opérateur de transfert} associé à $\theta\in\Ec^*$ est défini, 
pour $\gamma$ assez petit, et pour $\ab[\Re(\theta)]<\epsilon_0$ par 
\[\app[P_{\theta}]{\Hg}{\Hg}{\phi}
{\eta\mapsto\int_{g\in G}  e^{\theta\cdot\sgb(g,\,\eta)}\phi(g\cdot \eta)\dd \mu(g)}.\]
\end{df}
Pour $\theta=0$, on a le résultat suivant :
\begin{prop} Notons $N$ l'application linéaire  
\[\app[N]{\Hg}{\C}{\phi}{\int_{\flag}\phi(\eta)\dd\nu(\eta)}.\]
Alors l'opérateur $P_0$ vérifie les propriétés suivantes :
\begin{itemize}
\item  $\n[P_0]\leq 1$
\item  $P_0\fcar=\fcar$
\item $\n[P_{0\; |\Ker N}]<1$.
\end{itemize}

\end{prop}

Une preuve se trouve dans \cite[\S 10.9]{Livre}.
Cette proposition s'étend aux opérateurs $P_{\theta}$, 
pour $\theta$ dans un voisinage de $0$,
par le théorème d'analyse fonctionnelle suivant, 
dont une preuve est fournie dans \cite[\S 10.16]{Livre}, par exemple.

\begin{prop}\label{caroptrs} 
Il existe un réel $\gamma >1$, 
un réel $\epsilon >0$, 
un voisinage ouvert $U$ de $0$ dans $\Ec^*$, 
des fonctions 
\[\app[\Lambda]{U}{\C}{\theta}{\lambda_{\theta}},\]
\[\app[\Phi]{U}{\Hg}{\theta}{\phi_{\theta}},\]
analytiques sur $U$, 
et, pour tout $\theta$ dans $U$, une projection  
\[N_{\theta}:\Hg\rightarrow \C\phi_{\theta}\] 
commutant avec $P_{\theta}$ telles que 
\begin{itemize} 
\item $\lambda_0=1,\,\phi_0=\fcar,\, N_0=N$
\item $\forall \theta\in U,\, \optr\phi_{\theta}=\lambda_{\theta}\pth$
\item $\forall\theta\in U,\, \ab[\lth-1]\leq\epsilon$
\item $\forall\theta\in U,\, \nu(\pth)=1$
\item Le rayon spectral de l'endomorphisme ${\optr}_{|\Ker{N_{\theta}}}$ est inférieur à $1-\epsilon$. 
\end{itemize}
\end{prop}

On a, d'après \cite[\S 10.17]{Livre}, le développement suivant de $\lth$ :

\begin{prop}\label{devlambda}
La différentielle de $\Lambda$ en $0$ est $\smb$, 
et sa différentielle seconde en $0$ est $\Phi_{\mu}+\smb^2$ 
pour une certaine forme bilinéaire $\Phi_{\mu}$ définie positive sur $E^*$. 
\end{prop}
Enfin, le résultat suivant sur le rayon spectral des opérateurs $\optri$ est dû à Guivarc'h 
et Benoist-Quint, on en trouvera une preuve dans \cite{Livre} :

\begin{prop}\label{rayonspectraloptri}
Pour tout $\theta\in E^*\setminus\{0\}$, l'opérateur $\optri$ de $\Hg$ 
est de rayon spectral strictement inférieur à $1$.
\end{prop}

\subsection{Un second critère de transience}

On peut à présent montrer la proposition suivante :

\begin{prop}\label{trsd3sm} Si $\sm$ est non-nulle modulo $\Lie[a']$ 
(i.e. si le cocycle $\sgb$ est non-centré) 
ou si la dimension $d$ de $E$ est supérieure ou égale à $3$, 
alors la \maX est transiente.
\end{prop}

\begin{dem} On va distinguer les deux cas.

\subparagraph{Premier cas : $\sm\in\Lie[a']$}
Pour $k\in\N$, et $\eta\in \flag$, 
on note $\mukx$ la probabilité sur $E$ image de $\mu^{*k}$ par l'application $g\mapsto \sgb (g,\,\eta)$. 
On fixe un $\upsilon >0$ tel que la boule $V$ de centre $0$ et de rayon $\upsilon$ 
soit incluse dans l'ouvert $U$ donné par le théorème \ref{caroptrs}.
Pour montrer que la marche est transiente en tout point, il suffit de montrer que la quantité
\[F_{\eta}(C)=\sum_{k=0}^{\infty}\int_{g\in G}\fcar_{C}(\sgb(g,\,\eta))\dd\mu^{*k}(g)
=\sum_{k=0}^{\infty}\mukx(C)\]
est finie pour tout $\eta\in \flag$ et tout compact $C$ de $E$, d'après le lemme \ref{Gxtranscocycle}. 
Fixons un tel $\eta$ et un tel $C$. 
On identifie $E$ et $\R^d$.
Comme la fonction $\fcar_C$ est à support compact et valeurs positives, 
on peut trouver une fonction $\phi$ intégrable sur $\R^d$ majorant $\fcar_C$, 
de transformée de Fourier $\hat{\phi}$ positive à support compact . 
On obtient les inégalités :
\begin{align*}
F_{\eta}(C) &\leq \sum_{k=0}^{\infty}\int_E \phi(t)\dd\mukx (t)\\
& =  \sum_{k=0}^{\infty}\int_{E^*} \hat{\phi}(\theta)\widehat{\mukx}(\theta)\dd\theta. \\
\end{align*}
On peut calculer la transformée de Fourier $\widehat{\mukx}$ grâce aux opérateurs de transfert :
\[\widehat{\mukx}(\theta)=\int_{t\in\R^d}e^{i\theta\cdot t}\dd\mukx(t)
=\int_{g\in G}e^{i\theta\cdot \sgb(g,\,\eta)}\dd \mu^{*k} (g) = \optri^k 1 (x).\]
En divisant l'intégrale en deux parties, on obtient l'inégalité :
\[F_{\eta}(C)
\leq \sum_{k=0}^{\infty} \int_{\n[\theta] \geq \eta} \hat{\phi}(\theta)\widehat{\mukx}(\theta)\dd\theta 
+ \sum_{k=0}^{\infty} \int_{\n[\theta] < \eta} \hat{\phi}(\theta)\widehat{\mukx}(\theta)\dd\theta \]
En utilisant la proposition \ref{caroptrs}, 
on obtient, outre les fonctions $\Phi$ et $\Lambda$, deux fonctions analytiques sur $V$ : 
\[\app[\Psi]{V}{\Hg}{\theta}{\psi_{\theta}\in\Ker N_{\theta}}\]
et
\[\app[a]{V}{\C}{\theta}{a(\theta)}\]
 vérifiant, pour tout $\theta\in V$,
\begin{equation}
\label{decomp1}
\fcar=a(\theta)\pth +\psi_{\theta},
\end{equation}
où on a noté $\fcar$ la fonction constante égale à $1$. 
On obtient l'égalité, pour tout $\theta\in V\cap E$ et tout $k$ :
\[\optri^k1(\eta)= a(i\theta)\lthi^k\pthi(\eta)+\optri^k\psi_{i\theta}(\eta). \]
Notons 
\[I_1= \sum_{k=0}^{\infty} \int_{\n[\theta] \geq \upsilon, \theta\in E^*} 
\hat{\phi}(\theta)\widehat{\mukx}(\theta)\dd\theta,\]
\[I_2= \sum_{k=0}^{\infty} \int_{\n[\theta] < \upsilon, \theta\in E^*} 
\hat{\phi}(\theta)\optri^k\psi_{i\theta}(\eta)\dd\theta,\]
\[I_3= \sum_{k=0}^{\infty} \int_{\n[\theta] < \upsilon, \theta\in E^*} 
\hat{\phi}(\theta)a(i\theta)\lthi^k\pthi(\eta)\dd\theta.\]
On a
\[F_{\eta}(C)\leq \ab[I_1]+\ab[I_2]+\ab[I_3]\]

La première somme $\ab[I_1]$ est majorée 
par une série géométrique de raison strictement inférieure à $1$ :
en effet, d'après la proposition \ref{rayonspectraloptri}, 
le rayon spectral de $\optri$ pour tout $\n[\theta] \geq v$ est strictement inférieur à $1$,
et on a choisi $\hat{\phi}$ à support compact.
Cette première somme est donc finie. 

Etudions maintenant $\ab[I_2]$.
Comme l'application $\Psi$ est analytique, quitte à prendre une boule $V'$ plus petite,
on peut la supposer bornée. 
D'après la proposition \ref{caroptrs}, 
le terme $\optri^k\psi_{i\theta}(\eta)$ décroît exponentiellement quand $k$ augmente, 
et ce, uniformément en $\theta$. 
Comme $\hat{\phi}$ est bornée car à support compact, 
on en déduit que $\ab[I_2]$ est finie. 

Considérons à présent $\ab[I_3]$.
Les applications $a$ et $\Phi$ sont analytiques, 
et quitte à considérer une boule $V'$ plus petite, 
on peut supposer qu'elles sont bornées pour les normes appropriées.
Comme $\hat{\phi}$ est bornée, il existe une constante $A$ telle que : 
\begin{align*}
\ab[I_3]& \leq A \int_{\n[\theta]<\upsilon} \sum_{k=0}^{\infty} \ab[\lthi^k] \dd\theta \\
& \leq A \int_{\n[\theta]<\upsilon} \frac{1}{1-\ab[\lthi]}\dd\theta.
\end{align*}
D'après la proposition \ref{devlambda}, on a le développement suivant de la quantité $\lthi$,
pour $\theta\in E^*$ : 
\[\lthi=1+i\smb(\theta)-\frac{1}{2}(\phi_{\mu}(\theta)+\smb^2(\theta))+o(\n[\theta]^2).\]
Par hypothèse, la moyenne $\smb$ est nulle, et on obtient le développement de $1-\ab[\lthi]$,
pour $\theta\in E^*$ :
\[ 1-\ab[\lthi]=\frac{1}{2}\phi_{\mu}(\theta)+o(\n[\theta]^2).  \]
Comme $\phi_{\mu}$ est une forme quadratique définie positive sur $E^*$,
en choisissant $\upsilon$ assez petit,
l'intégrale $\int_{\n[\theta]<\upsilon} \frac{1}{1-\ab[\lthi]}\dd\theta$ converge pour $d\geq 3$, 
et donc $\ab[I_3]$ est finie. 
Ainsi, $F_{\eta}(C)$ est fini, ce qui conclut la preuve de la transience dans le cas $\sm\in\Lie[a']$.

\subparagraph{Second cas : $\sm\notin\Lie[a']$}   
Cette démonstration est semblable à celle du principe des grandes déviations, 
\cite[Lemme 8.5]{Livre}.
Pour $k\in\N$, et $\eta\in \flag$, notons $\stb(g,\eta)=\sgb (g,\,\eta)-k\smb$ pour tout $g\in G$, 
et notons cette fois $\mukx$ la mesure image de $\mu^{*k}$ 
par l'application $g\mapsto \stb (g,\,\eta)$. 
Soit $C$ un compact de $E$, soit $\eta\in \flag$, considérons la quantité $F_{\eta}(C)$ :
\begin{equation}
\label{fxksmb}
F_{\eta}(C)=\sum_{k=0}^{\infty}\int_{g\in G}\fcar_{C}(\sgb(g,\,\eta))\dd\mu^{*k}(g)
=\sum_{k=0}^{\infty}\mukx(C-k\smb).
\end{equation}
A nouveau, d'après le lemme \ref{Gxtranscocycle}, il suffit de montrer que $F_{\eta}(C)$ est finie.
Or, on a l'inclusion, pour tout $k\in\N$ assez grand, 
pour un certain $M>0$ ne dépendant que de $C$ et $\smb$ :
\[\{g\in G |\;\stb(g,\,\eta)\in C-k\smb\}\subset \{g\in G |\;\n[\stb(g,\,\eta)] \geq k M\}. \]
Notons $H_k= \{g\in G |\;\n[\stb(g,\,\eta)] \geq k M\}$ l'ensemble de droite. 
Choisissons un ensemble fini $\Theta\subset E^*$ 
(formé par exemple de $2d$ formes proportionnelles aux formes coordonnées) 
tel qu'on ait l'inclusion, pour tout $k\geq 1$ :
\begin{equation}\label{HkKtk}
H_k \subset \bigcup_{\theta\in\Theta} K_{\theta,\,k}, 
\end{equation}
en notant $K_{\theta,\,k}$ :
\[ K_{\theta,\,k} = \{g\in G |\;\theta(\stb(g,\,\eta)) \geq k \} .\]
Calculons la mesure $\mu^{*k}(K_{\theta,\,k} )$ pour un $\theta\in\Theta$.
L'inégalité de Bienaymé-Tchebychev donne, pour tout $t>0$:
\begin{equation}\label{bienaymetcheb} \mu^{*k}( K_{\theta,\,k})
=\mu^{*k}( \{g\in G |\; e^{t\theta(\stb(g,\,\eta))} \geq e^{t k} \})
\leq e^{-tk} \int_G e^{t\theta(\stb(g,\,\eta))}\dd\mu^{*k}(g). \end{equation}
Choisissons un voisinage $W$ de $0$ dans $\R$ vérifiant $W\theta\subset U$.  
L'inégalité (\ref{bienaymetcheb}) s'écrit encore, pour $t\in W$  :
\[\mu^{*k}( K_{\theta,\,k})\leq e^{-tk} e^{-tk\theta(\smb)}P_{t\theta}^k \fcar (\eta).\]
En écrivant la fonction constante $\fcar$ selon la décomposition (\ref{decomp1}), 
on obtient l'égalité :
\[P_{t\theta}^k\fcar(\eta)= a(t\theta)\lambda_{t\theta}^k\phi_{t\theta}(\eta)+
P_{t\theta}^k\psi_{t\theta}(\eta). \]
A nouveau, d'après la proposition \ref{caroptrs}, pour $t\in W$,
$\n[P_{t\theta}^k\psi_{t\theta}]$ décroît exponentiellement quand $k$ tend vers l'infini.
On en déduit une majoration de la limite supérieure de la quantité $ \log \n[P_{t\theta}^k \fcar(\eta)]$ :
\[\limsup_{k\tend\infty} \frac{1}{k} \log \n[P_{t\theta}^k\fcar(\eta)] \leq \log \lambda_{t\theta},  \]
ce qui donne la majoration asymptotique suivante de la quantité $\mu^{*k}( K_{\theta,\,k})$,
pour tout $t\in W$  :
\[\limsup_{k\tend\infty} \frac{1}{k} \log \mu^{*k}( K_{\theta,\,k}) 
\leq \log \lambda_{t\theta} - t(1+\theta(\smb)).\]
En $t=0$, la fonction
$\tau: t\mapsto  \log \lambda_{t\theta} - t(1+\theta(\smb))$ 
prend la valeur $0$ et sa dérivée est égale à $-1$. 
La fonction $\tau$ prend donc des valeurs strictement négatives en un point $t_0\in W\cap \R_+^*$. 
On en déduit que la limite $\limsup_{k\tend\infty} \frac{1}{k} \log \mu^{*k}( K_{\theta,\,k})$
est strictement négative. 
Notons-la $l_{\theta}$, choisissons $r_{\theta}$ tel que
 $e^{l_{\theta}}<r_{\theta}<1$.
On a alors la relation de domination :
\[\mu^{*k}(K_{\theta,\, k})=\underset{k\rightarrow \infty}{\mathcal{O}}(r_{\theta}^k) \] 
L'inégalité (\ref{HkKtk}) devient, pour tout $k\in\N$,
\[\mu^{*k}(H_k)\leq\sum_{\theta\in\Theta} \mu^{*k}(K_{\theta,\,k}).\]
Comme l'ensemble $\Theta$ est fini, on en déduit la domination : 
\[\mu^{*k}(H_k)=\underset{k\rightarrow \infty}{\mathcal{O}}(\sum_{\theta\in\Theta} r_{\theta}^k),\]
et donc 
\[ \mukx(K-k\smb)=\underset{k\rightarrow \infty}{\mathcal{O}}(\sum_{\theta\in\Theta} r_{\theta}^k).\]
La série de terme général $\mukx(K-k\smb)$ est donc dominée 
par une somme finie de séries géométriques convergentes.
La fonction de Green $F_{\eta}(C)$ est donc finie, ce qui conclut la preuve.
\end{dem}

\section{Une condition suffisante de récurrence uniforme}

Dans cette partie, nous allons montrer une condition suffisante de récurrence uniforme
dans la proposition \ref{d2rec}
qui permettra de conclure la preuve des théorèmes \ref{thmdicho} et \ref{thmclass}.

\subsection{Récurrence presque partout}

\subsubsection{Théorème de Schmidt-Conze}

Un théorème de Schmidt, démontré dans  \cite{Schm98}, 
qui améliore un résultat de Conze, dans \cite{Conz}, 
va permettre d'obtenir une condition de récurrence du cocycle d'Iwasawa.
On se place sur $(Z,\, \mathcal{Z},\, \lambda)$, un espace borélien probabilisé, 
muni d'un automorphisme $T:Z \rightarrow Z$ préservant la mesure, supposée ergodique.
On choisit un entier $d\geq 0$.
On considère alors une application borélienne $f:Z\rightarrow \R^d$, 
et on munit $\R^d$ d'une norme $\n[\cdot]$.
On définit les sommes de Birkhoff de $f$, pour $z\in Z$ et $n\in\N$ : 
\[f(n,z)=\left\{\begin{array}{rl} \sum_{k=0}^{n-1}f(T^k z) & \text{ si } n\geq 1\\
0 & \text{ si } n=0 \end{array}\right. .\]

\begin{df}\label{frec} On dit que $f$ est \emph{récurrente en presque tout point de $Z$} 
si on a :
\[\lambda(\{z\in Z \,|\, \liminf_{\nti} \n[f(n,z)]=0\})=1 . \]

On dit que $f$ \emph{vérifie le \tcl} si la suite de variables aléatoires 
$Z_n=\frac{f(n,\,\cdot)}{\sqrt{n}}$ converge en loi vers une loi normale, éventuellement dégénérée.
\end{df}

Énonçons à présent le théorème de Schmidt-Conze :

\begin{prop}\label{schmidt} (Schmidt-Conze) 
Soit $(Z,\, \mathcal{Z},\, \lambda)$ un espace borélien probabilisé, 
et $T:Z\rightarrow Z$ un automorphisme ergodique préservant la mesure. 
Supposons $d\leq 2$.
Soit $f: Z\rightarrow \R^d$ une application borélienne vérifiant le \tcl. 
Alors $f$ est récurrente en presque tout point de $Z$. 
\end{prop}

\subsubsection{Construction d'un système dynamique inversible} 

Nous allons nous ramener à la proposition \ref{schmidt} 
pour montrer une condition suffisante de récurrence presque partout de la \maX.

Notons $d$ la codimension de $A'$ dans $A$. 
On peut alors assimiler l'espace $E=\Lie[a]/\Lie[a']$ à $\R^d$. 
Il est naturel de vouloir considérer le système dynamique probabilisé  
$(\Bf,\,\borBf,\,T,\,\beta\otimes\nu)$,
où $\borBf$ est la tribu borélienne de $\Bf$, et
où $T$ est l'application préservant la mesure $\beta\otimes\nu$ :
\[\app[T]{\Bf}{\Bf}{(b,\,\eta)}{(Sb,\,b_1\eta)},\]
en notant
\[\app[S]{B}{B}{b=(b_1,\,b_2,\hdots)}{Sb=(b_2,\,b_3,\,\hdots)}\] 
le shift sur $B$.
Ce système dynamique est ergodique, puisque la mesure $\nu$ est $\mu$-ergodique.
Une construction due à Furstenberg (\cite{Furst63}) en fournit une extension inversible.
Notons $\Bb$ l'espace $\Bb=G^{\Z}$, $\mbb$ la mesure produit $\mu^{\otimes\Z}$, 
et $\borBb$ la tribu borélienne de $\Bb$. 
Notons $\tilde{S}$ le shift sur l'espace de Bernoulli bilatère $(\Bb,\,\borBb)$.
Rappelons la construction suivante de Furstenberg : 

\begin{lem}\label{nub} 
Pour $\mbb$-presque tout $b\in\Bb$, la suite $((b_1\cdots b_n)_*\nu)\sui$
de mesures de probabilité sur $\flag$ a une limite $\nu_b$.
Cette limite est la mesure de Dirac en un point que l'on notera $\xi_b$.
On a de plus l'égalité
\begin{equation}\label{nubeq}\nu=\int_{\Bb}\nu_b\dd\mbb(b).\end{equation}
\end{lem}

On en trouvera une démonstration dans \cite{Livre}. 
Pour tout $b\in\Bb$, on note $b_+$ la suite $(b_1,\,b_2,\hdots,\,b_n,\hdots)$,
et $b_-$ la suite $(\hdots,\,b_{-n},\hdots,\,b_{-1},\,b_0)$. 
Notons $\nu_{b_-}=\delta_{\xi_{b_-}}$, pour $\mbb$-presque tout $b\in\Bb$, 
la mesure de probabilité limite sur $\flag$ suivante :
\[\nu_{b_{-}}=\lim_{\nti}(b_0\cdots {b_{-n}})_{*}\nu.\]
En particulier, pour $\mbb$-presque tout $b\in\Bb$, on a 
$\xi_{(Sb)_-}=b_1\xi_{b_-}$. 
Le lemme suivant se déduit directement de ces considérations :

\begin{lem} 
L'application 
\[\app[\rho]{\Bb}{B\times\flag}{b}{(b_+,\,\xi_{b_-})}\]
est un morphisme de systèmes dynamiques au sens suivant :
on a $T\circ\rho = \rho\circ \tilde{S}$,
et l'image par $\rho$ de la mesure $\mbb$ est $\beta\otimes\nu$. 
\end{lem}

On se ramène donc à l'étude du système dynamique $(\Bb,\,\borBb,\,\tilde{S},\,\mbb)$,
dont $(\Bf,\,\borBf,\,T,\,\beta\otimes\nu)$ peut être vu comme un facteur via $\rho$. 

\subsubsection{Application du théorème de Schmidt-Conze}

Notons $\tilde{f}$ l'application borélienne suivante :
\[\app[\tilde{f}]{\Bb}{\R^d}{b}{\sgb(b_1,\,\xi_{b_-})}. \]
Puisque $\sgb$ est un cocycle, on constate immédiatement la propriété suivante, 
étant donnés $b\in\Bb$ et $n$ un entier strictement positif :
\[\tilde{f}(n,\,b) = \sum_{k=0}^{n-1}\tilde{f}(T^k(b,\,\eta)) = \sgb(b_n\cdots b_1,\, \xi_{b_-}). \]

Rappelons le résultat suivant dû à Goldsheid et Guivarc'h (voir \cite{Goldguiv}) 
qui constitue un \tcl pour le cocycle d'Iwasawa.

\begin{prop} (Théorème central-limite pour le cocycle d'Iwasawa)\label{tcliwas}
Avec les hypothèses et notations ci-dessus, 
supposons que la mesure $\mu$ ait un moment exponentiel fini. 
Alors il existe une loi normale de probabilité $N_{\mu}$ sur $E$ 
telle que pour toute fonction continue et bornée $\psi$ sur $E$ 
on ait la convergence, uniforme en $\eta\in\flag$ :
\[\int_{g\in G}\psi(\frac{\sgb (g,\,\eta)-n\smb}{\sqrt{n}})\dd\muen(g)
\xrightarrow[\nti]{}\int_{E}\psi\dd N_{\mu}.  \]
\end{prop}

Munis de ces résultats, nous pouvons démontrer le lemme suivant : 

\begin{lem}\label{d2recpartielle} Si $H$ est de la forme $A'N$, 
où $A'$ est un sous-groupe de $A$ de codimension au plus $2$, 
et si l'algèbre de Lie $\Lie[a']$ de $A'$ contient la moyenne $\sm$ du cocycle d'Iwasawa, 
alors pour $\beta\otimes\nu$-presque tout $(b,\,\eta)\in\Bf$ on a
\[\liminf_{\nti}\n[\sgb(b_n\cdots b_1,\,\eta)]=0.\] 
\end{lem}

\begin{dem} On applique le théorème de Schmidt-Conze (proposition \ref{schmidt}) 
au système dynamique $(\Bb,\,\mathcal{\Bb},\,\tilde{S},\,\beta)$, muni de l'application $\tilde{f}$,
qui vérifie le théorème central-limite d'après la proposition \ref{tcliwas}, 
puisque, par hypothèse, $\smb$ est nulle.
On obtient, pour $\mbb$-presque tout $b\in \Bb$ :
\[\liminf_{\nti} \n[\sgb(b_n\cdots b_1,\,\xi_{b_-})]=0, \]
et donc, pour $\beta\otimes\nu$-presque tout $(b,\,\eta)\in\Bf$ :
\[\liminf_{\nti} \n[\sgb(b_n\cdots b_1,\,\eta)]=0. \]
\end{dem}

\subsection{Un critère de récurrence uniforme}

Le lemme \ref{d2recpartielle} ne suffit pas encore à donner la récurrence sur tout $X$ de la \maX, car il ne donne cette information que
pour $\nu$-presque tout $\eta$ dans $\flag$.
Nous allons donc montrer la proposition suivante.

\begin{prop}\label{d2rec} Si $H$ est de la forme $A'N$, 
où $A'$ est un sous-groupe de $A$ de codimension au plus $2$, 
et si l'algèbre de Lie $\Lie[a']$ de $A'$ contient la moyenne $\sm$ du cocycle d'Iwasawa, 
alors la \maX est uniformément récurrente sur tout $X$. 
\end{prop}
 
D'après le lemme \ref{reccocycleiwas}, 
il suffit de montrer que pour tout $\eta\in\flag$, 
pour $\beta$-presque tout $b\in B$, pour tout $t\in E$, on a
\[\liminf_{\nti}\n[t+\sgb(b_n\cdots b_1,\,\eta)] = 0,\]
c'est à dire que la suite $(\sgb(b_n\cdots b_1,\,\eta))_{n\in\N}$
est dense dans $E$.  
Pour ce faire, on montrera que 
pour $\beta\otimes\nu$-presque tout $(b,\,\eta)\in B\times\flag$,\
la suite $(\sgb(b_n\cdots b_1,\,\eta))\sui$ est dense dans $\R^d$
(c'est le lemme \ref{systfibredense}) ; 
on montrera également que pour tous $\eta,\,\eta'\in\flag$, la suite
$(\sigma(b_n\cdots b_1,\,\eta)-\sigma(b_n\cdots b_1,\,\eta'))\sui$
a une limite finie (c'est le lemme \ref{iwakapdiffbornee}).
L'association de ces deux résultats permettra de conclure. 
Commençons par définir la projection de Cartan $\kappa$ sur $G$.

\begin{df} Écrivons la décomposition de Cartan du groupe $G$ : $G=K\exp \Lie[a]_+ K$. 
Pour tout $g\in G$, il existe un unique élément $\kappa(g)\in\Lie[a]_+$ tel qu'on ait
\[g\in K\exp(\kappa(g)) K.\]
L'application $\kappa : g\mapsto \kappa(g)$ est appelée la \emph{projection de Cartan}. 
\end{df} 

\begin{lem}\label{iwakapdiffbornee} Pour tout $\eta\in\flag$, 
pour $\beta$-presque tout $b\in B$,
la suite $(\sigma(b_n\cdots b_1,\,\eta)-\kappa(b_n\cdots b_1))\sui$
a une limite $\ell_{b,\,\eta}$. 
En particulier, pour tous $\eta,\,\eta'\in\flag$, la différence
$\sigma(b_n\cdots b_1,\,\eta)-\sigma(b_n\cdots b_1,\,\eta')$
tend vers la limite finie $\ell_{b,\,\eta}-\ell_{b,\,\eta'}$. 

\end{lem}

\begin{dem} Rappelons le résultat suivant (prouvé dans \cite{Livre}) : 
Soit $(V,\rho)$ une représentation irréductible de $G$, de plus haut poids $\chi$. 
Notons $V_{\chi}$ l'espace de poids associé à $\chi$, et pour $\eta=gP\in\flag$, 
notons $V_{\eta}$ le sous-espace $V_{\eta}=g\cdot V_{\chi}$ de $V$.
Alors il existe une norme euclidienne sur $V$, notée $\n[\cdot]$, vérifiant, 
pour tous $g\in G,\,\eta\in\flag,\, v\in V_{\eta}$ :
\begin{itemize}
\item $\chi(\kappa(g))=\log \n[\rho(g)]$
\item $\chi(\sigma(g,\,\eta))=\log\frac{\n[\rho(g)v]}{\n[v]}.$
\end{itemize}
Fixons une telle représentation, qu'on choisit de plus proximale. 
Montrons que la suite $(\chi(\sigma(b_n\cdots b_1,\,\eta)-\kappa(b_n\cdots b_1)))\sui$ 
a une limite pour $\beta$-presque tout $b\in B$,
pour tout $\eta\in\flag$.
Un lemme classique de Furstenberg (voir par exemple \cite[Prop 3.7]{Livre}) 
dit que pour $\beta$-presque tout $b\in B$, 
il existe un hyperplan $W_b\subset V$ tel que 
\begin{itemize}
\item pour toute valeur d'adhérence $\pi$ de la suite 
$(p_n=\frac{\rho(b_n\cdots b_1)}{\n[\rho(b_n\cdots b_1)]})\sui$, on a
$\Ker\pi=W_b$, 
\item pout tout $v\in V$ non-nul, on a $\beta(\{b\in B\, |\, v\in W_b \})=0.$
\end{itemize}
Soit $v\in V$ non-nul, soit $b\in B$ tel que $v\notin W_b$. 
Alors la quantité $\n[\pi(v)]$, strictement positive,
ne dépend pas de la valeur d'adhérence $\pi$
de la suite $(p_n)\sui$ choisie. 
En effet, fixons $f_b$ une forme linéaire non-nulle sur $\R^d$, de noyau $W_b$ et de norme $1$.
Les valeurs d'adhérence $\pi$
de la suite $(p_n)\sui$
s'écrivent alors $\pi\,:\,v'\mapsto f_b(v') v_{\pi}$ pour un certain vecteur $v_{\pi}$ de norme $1$,
puisque ces applications sont de rang $1$, de norme $1$ et de noyau $W_b$. 
La quantité $\n[\pi(v)]$ est donc toujours égale à $\ab[f_b(v)]$. 
La suite bornée $(\n[p_n(v)])\sui$ a une seule valeur d'adhérence, elle converge donc. 
La suite $(\log\frac{\n[\rho(b_n\cdots b_1)v]}{\n[v]}-\log \n[\rho(b_n\cdots b_1)])\sui$ 
converge alors.
Ceci étant vrai pour toute représentation proximale de $G$, 
on a montré le lemme. 
\end{dem}

\begin{lem}\label{systfibredense} Soit $d\leq 2$. 
Pour $\beta\otimes\nu$-presque tout $(b,\,\eta)\in B\times\flag$, 
la suite $(\sgb(b_n\cdots b_1,\,\eta))\sui$ est dense dans $\R^d$.
\end{lem}

La démonstration du lemme \ref{systfibredense} occupera la partie suivante.
On peut enfin prouver la proposition \ref{d2rec}. 

\paragraph{Preuve de la proposition \ref{d2rec}}
Choisissons grâce au lemme \ref{systfibredense}
un $\eta_0\in\flag$ tel que pour $\beta$-presque tout $b\in B$, 
la suite $(\sgb(b_n\cdots b_1,\,\eta_0))\sui$ est dense dans $\R^d$.
Alors, d'après le lemme \ref{iwakapdiffbornee}, 
pour $\beta$-presque tout $b\in B$, la suite $(\kappa(b_n\cdots b_1))\sui$ est dense dans $\R^d$. 
Et donc pour tout $\eta\in \flag$, pour $\beta$-presque tout $b\in B$, 
la suite $(\sgb(b_n\cdots b_1,\,\eta))\sui$ est dense.
Ainsi, pour $\beta$-presque tout $b\in B$, pour tout $t\in \R^d$, on a l'égalité :
\[\liminf_{\nti}\n[t+\sgb(b_n\cdots b_1,\,\eta)] = 0,\]
ce qui montre la proposition \ref{d2rec}. 
\qed 

\subsection{Preuve du lemme \ref{systfibredense}} 

On notera $\ell$ la mesure de Lebesgue sur $\R^d$. 
Considérons le système dynamique $\tilde{\Sigma}$ :
\[\tilde{\Sigma}=(\Bb\times\R^d,\,R,\, \mbb\otimes\ell),\]
avec $\Bb=B^{\Z}$, $\mbb=\beta^{\otimes\Z}$, $\tilde{S}$ le shift bilatère sur $\Bb$,
et $R$ l'application 
\[\app[R]{\Bb\times\R^d}{\Bb\times\R^d}{(b,\,t)}{(Sb,\,t+\sgb(b_1,\,\xi_{b_-}))}.\]
L'application $R$ préserve la mesure $\mbb\otimes\ell$.
Pour montrer le lemme \ref{systfibredense}, 
on commencera par remarquer que 
le système dynamique $\tilde{\Sigma}$ est conservatif et ergodique :
c'est l'objet de la proposition \ref{tildeSigmaconserg}.
Pour montrer l'ergodicité de $\tilde{\Sigma}$,
on s'intéressera, dans la partie \ref{ergauxiliaire},
à celle d'un système dynamique auxiliaire, 
avec des méthodes calquant celles de Guivarc'h (\cite{Guiv89})
Ces deux propriétés de conservativité et d'ergodicité impliquent, 
d'après la proposition \ref{ergimpliquedense}, 
la densité des trajectoires dans l'espace $\Bb\times\R^d$. 
Ce dernier résultat permettra de conclure la preuve du lemme 
\ref{systfibredense}
à la fin de la partie \ref{tildesigmadense}.  

\begin{df} Soit $(X,\,T,\,m)$ un sytème dynamique, avec $X$ un espace topologique,
$m$ une mesure borélienne non-nécessairement finie,
$T$ une application borélienne préservant la mesure. 
Le système $(X,\,T,\,m)$ est dit \emph{conservatif} si pour tout borélien $A\subset X$
de mesure $m(A)>0$, il existe $n\geq 0$ tel que 
\[m(T^{-n}A\cap A)>0. \]
Le système $(X,\,T,\,m)$ est dit \emph{ergodique} si pour tout borélien $A\subset X$ tel que
$T^{-1}A=A$, on ait soit $m(A)=0$, soit $m(A^c)=0$,
où $A^c$ dénote le complémentaire de $A$ dans $X$. 
\end{df}

\begin{prop}\label{tildeSigmaconserg}
Supposons $d\leq 2$. Le système dynamique $\tilde{\Sigma}$ est conservatif et ergodique.
\end{prop}

\begin{dem}
On va se ramener, via le lemme \ref{ergconseq},
à l'étude de la conservativité de $\tilde{\Sigma}$ 
et de l'ergodicité de $\check{\Sigma}$,
en notant $\check{\Sigma}$ un système 
dont l'extension naturelle universelle est $\tilde{\Sigma}^{-1}$. 
On pose 
\[\check{\Sigma}=(B \times\R^d,\,V,\,\beta\otimes\ell).\]
La mesure produit $\beta\otimes\ell$ préserve l'application
\[\app[V]{B\times\R^d}{B\times\R^d}
{(b,\,t)}{(Sb,\,t+\sgb(b_1^{-1},\,\xi_b))}.\]
Par la propriété de cocycle de $\sgb$, pour $b\in B$, $\eta\in\flag$,
la quantité $t+\sgb(b_1^{-1},\,\eta)$ s'écrit aussi $t-\sgb(b_1,\, b_1^{-1}\eta)$ ;
l'extension naturelle inversible de $\check{\Sigma}$ est donc bien
le système $\tilde{\Sigma}^{-1}=(\Bb\times\R^d,\,R^{-1},\,\mbb\otimes\ell).$

\begin{lem}\label{ergconseq} 
Si un des systèmes 
$\tilde{\Sigma},\,\check{\Sigma}$ est conservatif, alors les deux le sont.
Si un des systèmes 
$\tilde{\Sigma},\,\check{\Sigma}$ est conservatif et ergodique, 
alors les deux le sont.
\end{lem}

\begin{dem} Voir \cite[3.1]{Aaron97}. 
Le système $\tilde{\Sigma}^{-1}$ 
est l'extension naturelle inversible de $\check{\Sigma}$. 
La conservativité (resp. les conservativité et ergodicité) de $\tilde{\Sigma}$ 
et celle (resp. celles) de $\tilde{\Sigma}^{-1}$ 
sont équivalentes puisque $\tilde{\Sigma}^{-1}$ est l'inverse de $\tilde{\Sigma}$.
La conservativité (resp. les conservativité et ergodicité)  d'un système dynamique 
et celle (resp. celles) de son extension 
naturelle inversible sont équivalentes (\cite[Thm. 3.1.5]{Aaron97}).
\end{dem}

Reprenons la preuve de la proposition \ref{tildeSigmaconserg}. 
Le lemme \ref{d2recpartielle} prouve la conservativité de $\tilde{\Sigma}$
pour $d\leq 2$
(cf. par exemple \cite[Thm 11.1]{Schm76}),
donc celle de $\check{\Sigma}$. 

\begin{lem}\label{ergBRdd2} 
Le système $\check{\Sigma}$ est ergodique si $d=1$ ou $d=2$. 
\end{lem} 

La preuve du lemme \ref{ergBRdd2} fait l'objet de la partie \ref{ergauxiliaire}.
Pour $d\leq 2$, le système $\tilde{\Sigma}$ est donc conservatif et ergodique. 
\end{dem}

\subsubsection{Ergodicité du système dynamique $(B\times\R^d,\,\beta\otimes\ell,\, V)$}
\label{ergauxiliaire}
Cette sous-partie est consacrée à la preuve du lemme \ref{ergBRdd2}.
Ce lemme est presque une conséquence du résultat de Guivarc'h \cite[\S 2.3, Cor. 3]{Guiv89},
à quelques conditions près sur l'application $V$. 
Nous redonnons ici les principaux éléments de sa démonstration,
adaptés à notre cas. 
On utilise la caractérisation suivante de l'ergodicité de $(B\times\R^d,\,\beta\otimes\ell,\, V)$, 
dont la preuve sera rappelée en annexe (Proposition \ref{conderggen}) :

\begin{lem}\label{condexrel} (Guivarc'h) 
Notons $\Lun_0(\ell)=\{\phi\in\Lun(\ell)\,|\,\int_{\R^d}\phi\dd\ell=0\}$.
Pour toute application $\alpha=\psi\otimes\zeta$, 
avec $\psi\in\Lun(\beta)$
et $\zeta\in\Lun_0(\ell)$,
on définit la suite des poussés en avant $V^n_{*}\alpha$ de $\alpha$ par $V^n$
en notant, pour $n\in\N^*$, pour $(b,\,t)\in B\times\R^d$ :
\[V^n_{*}\alpha (b,\,t)=
\int_{G^n} \psi((g_n,\hdots,\,g_1,\,b_1,\,b_2,\hdots))\zeta(t+\sgb(g_n\cdots g_1,\,\xi_b))
\dd \mu^{*n}(g_1,\hdots,\, g_n).\]
Si pour toute telle fonction $\alpha$
la suite des poussés en avant de $(V^n_{*}\alpha)\sui$ 
a une limite nulle en norme $\Lun$,
\cad
\[\lim_{\nti}\nLun[V^n_{*}(\psi\otimes\zeta)]=0,\]
alors le système $(B\times\R^d,\,\beta\otimes\ell,\, V)$ est ergodique. 
\end{lem}

On en déduit le lemme \ref{ergBRdd2}.
Soient $\psi\in\Lun(\beta)$
et $\zeta\in\Lun_0(\ell)$.
D'après le lemme \ref{condexrel}, 
il suffit de prouver qu'on a la limite dans $\Lun$ : 
\begin{equation}\label{majnormeL1}
\nLun[{V}^n_{*}(\psi\otimes\zeta)]\xrightarrow{\nti}0.
\end{equation}
Définissons, pour $\beta$-presque tout couple $b\in B$,
la mesure borélienne $\mu_{n,\,\psi}^b$ sur $\R^d$,
qui a une fonction $h$ mesurable sur $\R^d$ associe la valeur 
\[\mu_{n,\,\psi}^b(h)= \int_{G^n}h(-\sgb(g_n\cdots g_1,\,\xi_b))
\psi((g_n,\hdots,\, g_1,\,b_1,\,b_2,\hdots))
\dd\mu^{\otimes n}(g_1,\hdots,\,g_n).
\]
Ecrivons, pour $f:B\rightarrow\R$ et $h:\R^d\rightarrow\R$ 
deux fonctions mesurables bornées,
la quantité 
$(V_*(\psi\otimes\zeta))(f\otimes h)=\int_{B\times\R^d}f(b)h(t)
V_{*}(\psi\otimes\zeta) (b,\,t)\dd\beta\otimes\ell(b,\,t)$ :
\begin{align*}
(V_*(\psi\otimes\zeta))(f\otimes h)&=\int_{B\times\R^d\times G}
f(b)h(t)
\psi((g,\,b_1,\,b_2,\hdots)) 
\zeta(t+\sgb(g,\,\xi_b))
\dd\mu\otimes\beta\otimes\ell(g,\,b,\,t)\\
&=\int_{B\times\R^d}
f(b)h(t)
(\mu_{1,\,\psi}^b*\zeta(t))\dd\beta\otimes\ell\\
\end{align*}
On obtient donc, pour tout $n\in\N$, 
l'égalité :
\[
(V^n_*(\psi\otimes\zeta))(f\otimes h)=\int_{B\times\R^d}
f(b)h(t)
(\mu_{n,\,\psi}^b*\zeta(t))\dd\beta\otimes\ell \]
et, pour $\beta\otimes\ell$-presque tout $(b,\,t)\in \Bb\times\R^d$, 
\[(V^n_*(\psi\otimes\zeta))(b,\,t)=(\mu_{n,\,\psi}^b*\zeta)(t)\]
On obtient une majoration en norme $\Lun$ :
\[\nLun[V^n_*(\psi\otimes\zeta)]\leq
\int_{B}\n[\mu_{n,\,\psi}^b*\zeta]_{\Lun(\ell)}\dd\beta. \]
Par convergence dominée, pour montrer la convergence (\ref{majnormeL1}),
il suffit de montrer, pour $\beta$-presque tout $b\in B$, la convergence 
\begin{equation}\label{majnormeL1bis}
\n[\mu_{n,\,\psi}^b*\zeta]_{\Lun(\ell)}\xrightarrow{\nti}0.
\end{equation}

Montrons (\ref{majnormeL1bis}) pour des fonctions 
$\psi\in\Lun(\beta)$ 
et $\zeta\in\Lun_0(\ell)$ vérifiant les propriétés suivantes :
\begin{itemize}
\item $\psi$ est $\gamma$-hölderienne pour un $\gamma>0$ ;
\item $\psi$ est à support compact ;
\item $\zeta$ est à décroissance exponentielle ;
\item la transformée de Fourier $\hat{\zeta}$ de $\zeta$ est à support compact ;
\item $\hat{\zeta}$ est nulle sur un voisinage de $0$.
\end{itemize}
Le sous-espace vectoriel de $\Lun(\beta)$ 
engendré par de telles applications $\psi$ 
est dense ;
de même, les applications $\zeta$ de la forme ci-dessus sont denses dans $\Lun_0(\ell)$.
On aura donc montré la convergence (\ref{majnormeL1bis}) pour $\psi$ et $\zeta$ quelconques.
Fixons $\psi$ et $\zeta$ avec les propriétés indiquées ci-dessus. 
Notons $C$ le support (compact) de $\psi$. 
Il existe un réel $R_0>0$ tel que l'image par $r\,:\,b\mapsto\sgb(b_1^{-1},\,\xi_b)$
de $C$ est incluse dans la boule de rayon $R_0$. 
Fixons $\epsilon>0$.
Calculons :
\begin{align*}
\n[\mu_{n,\,\psi}^b*\zeta]_{\Lun(\ell)}&=\int_{\R^d}\ab[\mu_{n,\,\psi}^b*\zeta(t)]\dd\ell(t)\\
&\leq \int_{\overline{B}(n(R_0+\epsilon))}\ab[\mu_{n,\,\psi}^b*\zeta(t)]\dd\ell(t) + 
\int_{\R^d\setminus\overline{B}(n(R_0+\epsilon))}\ab[\mu_{n,\,\psi}^b*\zeta(t)]\dd\ell(t)\\
&=I_1+I_2.
\end{align*}
L'inégalité de Schwarz appliquée à l'intégrale $I_1$ donne :
\begin{align*}
I_1=\int_{\overline{B}(n(R_0+\epsilon))}\ab[\mu_{n,\,\psi}^b*\zeta(t)]\dd\ell(t)&
\leq 2n(R_0+\epsilon)\int_{\overline{B}(n(R_0+\epsilon))}\ab[\mu_{n,\,\psi}^b*\zeta(t)]^2\dd\ell(t) \\
&\leq 2n(R_0+\epsilon) \int_{(\R^d)^*}\widehat{\mu_{n,\,\psi}^b}(s)\hat{\zeta}(s)\dd s.
\end{align*}
Pour un $s\in(\R^d)^*$, la transformée de Fourier de $\mu_{n,\,\psi}^b$ prend pour valeur :
\begin{align*}
\widehat{\mu_{n,\,\psi}^b}(s) &=\int_{\R^d} e^{-i\,s(t)}\dd\mu_{n,\,\psi}^b(t)\\
&=\int_{G^n}e^{i\,s(\sgb(g_n\cdots g_1,\,\xi_b))}
\psi((g_n,\hdots,\, g_1,\,b_0,\,b_1,\hdots))
\dd\mu^{\otimes n}(g_1,\hdots,\,g_n).\\
\ab[\widehat{\mu_{n,\,\psi}^b}(s)]&\leq 
\ninf[\psi]\ab[ \int_Ge^{i\,s(\sgb(g,\,\xi_b))}\dd\mu^{*n}(g)]\\
&\leq \ninf[\psi] \ab[P^n_{is}(\fcar)(\xi_b)],
\end{align*}
où on a noté $P_{is}$ l'opérateur de transfert associé à $is$,
et $\fcar$ la fonction constante égale à $1$. 
D'après les propositions \ref{caroptrs} et \ref{rayonspectraloptri}, 
cet opérateur est de rayon spectral strictement inférieur à $1$, 
et variant continûment avec $s$, pour tout $s$ non nul. 
Comme le support $\Supp \hat{\zeta}$ de $\hat{\zeta}$ est un compact ne contenant pas $0$, 
il existe un $a,\,0<a<1$, une constante $c(\psi)>0$, tels que pour tout $s\in\Supp\hat{\zeta}$, 
on ait l'inégalité, en notant $\ab[\cdot]_{\gamma}$ la norme Hölder sur $\mathcal{H}^{\gamma}$,
\begin{equation*}
\ab[\widehat{\mu_{n,\,\psi}^b}(s)]\leq
\ninf[\psi] \ninf[P^n_{is}(\fcar)]\leq
\ninf[\psi] \ab[P^n_{is}(\fcar)]_{\gamma}\leq
c(\psi)a^n.
\end{equation*}
On en déduit une majoration de $I_1$ :
\begin{equation}\label{majI1}
I_1\leq 2n(R_0+\epsilon)\ninf[\hat{\zeta}] c(\psi) a^n.
\end{equation}
Considérons à présent l'intégrale $I_2$. 
Puisque $\zeta$ est à décroissance exponentielle, 
il existe une constante $A>0$ telle que, pour tout $t\in\R^d$, on ait la majoration :
\[\zeta(t)\leq A e^{-\ab[t]}.\]
Calculons :
\begin{align*}
I_2&=\int_{\ab[t]\geq n(R_0+\epsilon)}\ab[\mu_{n,\,\psi}^b*\zeta(t)]\dd t\\
&=\int_{\ab[t]\geq n(R_0+\epsilon)}\int_{s\in\R^d} \ab[\mu_{n,\,\psi}^b(s)\zeta(t-s)]\dd t \dd s\\
&\leq A \ninf[\psi] \int_{\ab[t]\geq n(R_0+\epsilon)}\int_{\ab[s]\leq n R_0} e^{-\ab[t-s]} \dd t \dd s\\
&\leq 2d A \ninf[\psi] \int_{ t \geq n(R_0+\epsilon)}\int_{s \leq n R_0} e^{-t} e^s \dd t \dd s.
\end{align*}
On obtient enfin la majoration :
\begin{equation}\label{majI2}
I_2\leq 2d A \ninf[\psi]  e^{-n\epsilon}.
\end{equation}
En combinant les inégalités (\ref{majI1}) et (\ref{majI2}), on obtient, pour $n\in\N$ :
\[\n[\mu_{n,\,\psi}^b*\zeta]_{\Lun(\ell)}\leq 
2n(R_0+\epsilon)\ninf[\hat{\zeta}] c(\psi) a^n + 2d A \ninf[\psi]  e^{-n\epsilon}, \]
ce qui donne bien la limite (\ref{majnormeL1bis}). 
\qed
\bigskip

\subsubsection{Densité dans les fibres de $\tilde{\Sigma}$}\label{tildesigmadense}

Revenons à la preuve du lemme \ref{systfibredense}.
Rappelons la proposition très générale suivante, qui, dans notre cas particulier, 
donnera la densité des trajectoires dans les fibres. 

\begin{prop}\label{ergimpliquedense} 
Soit $X$ un espace localement compact, à base dénombrable d'ouverts, 
soit $m$ une mesure de Radon sur $X$ ,  
et soit $T:X\rightarrow X$ un endomorphisme inversible préservant $m$,
ergodique, et conservatif. 
Alors pour tout ouvert $U$ de $X$ de mesure non-nulle, pour $m$-presque tout $x\in X$,
pour tout $n_0\in \N$, il existe un $n\geq n_0$ tel que $T^n x\in U$. 
\end{prop}

\begin{dem} Soit $U$ un ouvert de $X$, $m(U)>0$. 
Notons $F_U$ l'ensemble des $x\in X$ revenant infiniment souvent dans $U$ sous l'action de $T$ :
\[F_U=\{x\in X\,|\, \exists (n_k)_{k\in\N}\in\Z^{\N}\,:
\, \forall k\in\N,\, n_k<n_{k+1}\;\text{et}\;T^{n_k}x\in U\}.\]
Cet ensemble étant $T$-invariant, et $T$ étant ergodique, on a l'alternative
\begin{itemize}
\item $m(F_U)=0$,
\item $m(X\setminus F_U)=0$.
\end{itemize}
Quitte à considérer les intersections de $U$ avec une suite exhaustive de compacts de $X$, 
on peut supposer $U$ de mesure finie. 
Alors la fonction $f=\fcar_U$ est intégrable et positive.
Comme $T$ est conservatif, on a, $m$-presque sûrement :
\[\sum_{k=0}^{\infty}T^k\fcar_U=\infty,\]
et donc on ne peut pas avoir $m(F_U)=0$,
ce qui prouve le résultat.
\end{dem}

\paragraph{Preuve du lemme \ref{systfibredense}}
En appliquant la proposition \ref{ergimpliquedense} au système $\tilde{\Sigma}$,
on obtient le lemme \ref{systfibredense}. 
Le système $\tilde{\Sigma}$ est conservatif et ergodique 
d'après la proposition \ref{tildeSigmaconserg} ; 
on en déduit que pour presque tout $(b,\,t)\in\Bb\times\R^d$, 
la suite des itérées $(R^n(b,\,t))\sui$ est dense dans $\Bb\times\R^d$.
On en déduit que pour $\mbb$-presque tout $b\in\Bb$,
la suite $(\sgb(b_n\cdots b_1,\,\xi_{b-})\sui$
est dense dans $\R^d$; 
pour $\beta\otimes\nu$-presque tout $(b,\,\eta)\in\Bf$,
la suite $(\sgb(b_n\cdots b_1,\,\eta)\sui$
est dense dans $\R^d$.
\qed

\section{Conclusion}\label{preuveprinc}

Nous avons à présent tous les éléments nécessaires à la démonstration 
du théorème de dichotomie \ref{thmdicho} et du critère de récurrence \ref{thmclass}.

\paragraph{Preuve des thèoremes \ref{thmdicho} et \ref{thmclass}}
Le théorème \ref{thmdicho} est un corollaire du théorème \ref{thmclass}. 
Démontrons ce dernier.
On peut supposer $H$ de la forme $A'N'$, 
avec $A'$ et $N'$ des sous-groupes de Lie algébriques de $A$ et $N$, 
d'après le lemme \ref{iwasH}.
On a alors les possibilités suivantes :

\subparagraph{Premier cas : $N'\subsetneq N$} 
On est dans le cadre de la proposition \ref{roledeN} : la \maX est transiente.

\subparagraph{Second cas : $N'=N$ et $\sm\notin\Lie[a']$} 
On est dans le cadre de la proposition \ref{trsd3sm} : la \maX est transiente.

\subparagraph{Troisième cas : $N'=N$, $\sm\in\Lie[a']$, 
et $A'$ est de codimension au moins $3$ dans $A$} 
On est dans le cadre de la proposition \ref{trsd3sm} : la \maX est transiente.

\subparagraph{Dernier cas : $N'=N$, $\sm\in\Lie[a']$, 
et $A'$ est de codimension au plus $2$ dans $A$} 
On est dans le cadre de la proposition \ref{d2rec} : la \maX est uniformément récurrente. 
\qed

\appendix

\section{Un résultat de Guivarc'h sur l'ergodicité}

Nous reprenons dans cette annexe des résultats dus à Guivarc'h (\cite{Guiv89}). 
L'objet est de donner une preuve du lemme \ref{condexrel} :
il se déduit immédiatement de la proposition \ref{conderggen} ci-dessous. 
Rappelons quelques définitions.
Soit $X$ un espace topologique muni de sa tribu borélienne et d'une mesure de probabilité $m$,
$T:X\rightarrow X$ un endomorphisme préservant la mesure de $X$, 
et $f:X\rightarrow \R^d$ une application mesurable, 
avec $\R^d$ muni de la mesure de Lebesgue $\ell$. 
On définit alors $E=X\times\R^d$ l'espace fibré correspondant,
muni de la mesure borélienne $m\otimes\ell$, 
et de l'endomorphisme préservant la mesure
\[\app[T_f]{E}{E}{(x,\,t)}{(Tx,\,t+f(x))}.\]
Notons $\mathcal{E}$ la tribu borélienne de $E$, 
et $v$ la tribu engendrée par les ensembles de la forme 
$U\times\R^d$, avec $U$ un ouvert de $X$.
On appelle $v$ la \emph{partition en verticales} de $E$.

\begin{df}\label{defexexrel} On dit qu'une fonction $F:E\rightarrow \R$ est \emph{asymptotique}
si pour tout $n\in \N$, il existe une fonction $F_n :E\rightarrow \R$ telle qu'on ait 
$m\otimes\ell$-presque partout l'égalité
$F=F_n\circ T_f^n$.
On dit que le système $(E,\,T_f,\,m\otimes\ell)$ est \emph{exact}
si la tribu queue, donnée par $\bigcap_{n\geq 0}T_f^n\mathcal{E}$, est triviale, 
\cad si toutes les fonctions asymptotiques sur $E$
sont constantes $m\otimes\ell$-presque partout. 
On dit que le système $(E,\,T_f,\,m\otimes\ell)$ est \emph{exact relativement à $v$} si 
toute fonction $F:E\rightarrow \R$ asymptotique 
ne dépend que de la première variable, 
\cad si pour $m$-presque tout $x\in X$,
pour $\ell$-presque tous $t,\,u\in\R^d$, on a l'égalité $F(x,\,t)=F(x,\,u)=F(x)$. 
\end{df}

Le lemme suivant découle immédiatement de ses définitions :

\begin{lem}\label{ergexrelerg} Si l'application $T$ est ergodique sur $X$, 
et que l'application $T_f$ est exacte relativement à $v$, 
alors l'application $T_f$ est ergodique sur $E$.
\end{lem}

Notons $\Lun_0(m\otimes\ell)$ l'ensemble des applications $\Lun$ d'intégrale nulle sur $E$, 
$\Lun_0(\ell)$ l'ensemble des applications $\Lun$ d'intégrale nulle sur $\R^d$,
et $\Lun_v(m\otimes\ell)$ l'ensemble des applications $\Lun$ d'intégrale nulle dans chaque fibre,
\cad
\[\alpha\in\Lun_v\Leftrightarrow \mathrm{Pour}\,m\,\mathrm{-presque}\, 
\mathrm{tout}\, x\in X,\, \int_{\R^d} \alpha(x,\,t) \dd t =0. \]
Pour toute application $\alpha\in\Lun(E)$, on définit son poussé en avant par $T_f$ : 
c'est l'application ${T_f}{_*}\alpha : E\rightarrow \R $ 
vérifiant 
\[\forall \gamma\in\Linf(E),\, \int_E {T_f}{_*}\alpha (z) \gamma(z) \dd (m\otimes\ell)(z) =
\int_E \alpha(z) \gamma (T_f (z))\dd (m\otimes\ell)(z).\]

On va montrer, grâce à une caractérisation élégante 
de l'exactitude relative proposée par Guivarc'h, la proposition suivante, 
dont se déduit immédiatement le lemme \ref{condexrel}.

\begin{prop}\label{conderggen} Supposons le système dynamique $(X,\,T,\,m)$ ergodique.
Si pour toute application $\alpha=\psi\otimes\xi$, 
avec $\psi\in\Lun(m)$
et $\xi\in\Lun_0(\ell)$,
la suite des poussés en avant de $\alpha$ par $T_f^n$ a une limite nulle en norme $\Lun$,
\cad
\[\lim_{\nti}\nLun[{T_f}^n{_*}(\psi\otimes\xi)]=0,\] 
alors le système $(E,\,T_f,\,m\otimes\ell)$ est ergodique. 
\end{prop}

\begin{dem} Cette proposition se déduit immédiatement des lemmes \ref{ergexrelerg},
\ref{condexrel1}, et \ref{condexrel2}.
\end{dem}

\begin{lem}\label{condexrel1}Si pour toute application $\alpha\in\Lun_v(m\otimes\ell)$, 
la suite des poussés en avant de $\alpha$ par $T_f^n$ a une limite nulle en norme $\Lun$,
\cad
\[\lim_{\nti}\nLun[{T_f}^n{_*}\alpha]=0,\]
alors le système $(E,\,T_f,\,m\otimes\ell)$ est exact relativement à $v$. 
\end{lem}

\begin{dem}(cf. \cite{Guiv89}) 
Soient $\alpha\in\Lun_v(m\otimes\ell)$ 
et $F$ bornée, asymptotique sur $E$. 
Notons $F_n$ les applications telles que pour tout $n\in\N$, on ait $F=F_n\circ T_f^n$.
Calculons :
\[\ab[\int_{E}F(x,\,t)\alpha(x,\,t)\dd m(x)\dd t]=
\ab[\int_E F_n(x,\,t) {T_{f}^n}{_*}\alpha(x,\,t)\dd m(x)\dd t ] \leq
\ninf[F]\nLun[{T_f^n}{_*}\alpha].\]
Par hypothèse, on obtient
\[\forall \alpha\in\Lun_v(m\otimes\ell),\,\int_{E}F(x,\,t)\alpha(x,\,t)\dd m(x)\dd t=0.\]
On en déduit, pour tout $\xi\in\Lun_0(\ell)$, pour tout $\psi\in\Lun(m)$, l'égalité
\[\int_{E}F(x,\,t)\psi(x)\xi(t)\dd m(x)\dd t=0,\]
et on en déduit, en choisissant les fonctions $\psi$ et $\xi$ adéquates, 
que $F$ ne dépend pas de $t$, ce qui montre l'exactitude relativement à $v$. 
\end{dem}

\begin{lem}\label{condexrel2} Les applications $\alpha$ de la forme $\alpha=\psi\otimes\xi$, 
avec $\psi\in\Lun(m)$
et $\xi\in\Lun_0(\ell)$,
forment une partie totale de $\Lun_v(m\otimes\ell)$.
\end{lem}

\begin{dem} Considérons la partie $A$ de $\Linf(m\otimes\ell)$ suivante :
\[A=\{g\in\Linf(m\otimes\ell)\,|\, \mathrm{Pour}\,m\,\mathrm{-p.t.}\, 
\, x\in X,\,l\,\mathrm{-p.t.}\, 
t,\, u \in \R^d,\, g(x,\,t)=g(x,\,u)=g(x) \}.\]
Montrons que $\Lun_v(m\otimes\ell)$ est l'orthogonal de $A$. 
Soit $\alpha\in\Lun_v(m\otimes\ell)$, soit $g\in A$.
Calculons leur produit scalaire :
\begin{align*}
\int_E\alpha g \dd m\otimes\ell  &=\int_{X}\int_{\R^d} g(x, \,t) \alpha(x,\,t) \dd t \dd m(x)\\
&=\int_X g(x) \int_{\R^d} \alpha(x,\,t) \dd t \dd m(x) \\
&=\int_X g(x) \cdot 0 \dd m(x)\\
&= 0.\\
\end{align*}
On a donc l'inclusion :
\[\Lun_v(m\otimes\ell)\subset A^{\perp}.\]
Considérons à présent $\alpha\in A^{\perp}$. 
En choisissant les $g$ de la forme $\fcar_{U\times\R^d}$, 
avec $U$ une partie borélienne de $X$, 
on obtient, pour $m$-presque tout $x\in X$,
l'égalité :
\[0=\int_{\R^d} \alpha(x,\,t) \dd t ,\]
et donc $\alpha$ est bien un élément de $\Lun_v(m\otimes\ell)$. 
Notons $C$ le sous-espace vectoriel de $\Lun_v(m\otimes\ell)=A^{\perp}$
engendré par les applications de la forme $\alpha=\psi\otimes\xi$, 
avec $\psi\in\Lun(m)$
et $\xi\in\Lun_0(\ell)$.
Soit $g\in C^{\perp}\subset\Linf(m\otimes\ell)$. 
Alors on a
\[\forall\psi\in\Lun(m),\,\forall\xi\in\Lun_0(\ell),\,
\int_E g(x,\,t)\psi(x)\xi(t) \dd m\otimes\ell(x,\, t) = 0.\]
En choisissant $\xi$ comme combinaison linéaire de fonctions caractéristiques, 
on obtient, 
pour $\ell$-presque tout $t\in\R^d$ :
\[\forall\psi\in\Lun(m),\,
\int_X g(x,\,t)\psi(x) \dd m(x) = c(\psi),\]
avec $c(\psi)$ une constante ne dépendant pas de $t$.
En choisissant pour $\psi$ des fonctions caractéristiques, 
on obtient alors que $g$ ne dépend pas de la variable $t$ pour $m$-presque tout $x\in X$, 
\cad $g\in A$. 
On obtient donc l'inclusion $C^{\perp}\subset A$,
l'inclusion réciproque étant immédiate. 
On obtient enfin la densité de $C$ dans $\Lun_v(m\otimes\ell)$.
\end{dem}

\nocite{Schott}

\bibliography{recbib}{}
\bibliographystyle{abbrv}

\end{document}